# FROM ART TO GEOMETRY: AESTHETIC AND BEAUTY IN THE LEARNING PROCESS

**PASTENA Nicolina, (I), PALLADINO Nicla, (I), VACCARO Maria Alessandra (I)**

**Abstract.** Starting from the concept that knowledge comes as element of mediation between the convergent thinking, founded on experience, and the divergent thinking, placed in the perceptive, intuitive, creative dimension, in this paper we want to present an idea for developing an educational path combining the concept of "beauty" and some historical notes. It is possible to use this dissertation as a starting point to conceive a geometric laboratory that drawing inspiration from artistic works, get to create geometric shapes provided with fascinating symmetries.

## 1   Creative thought and mathematical thought: rationality and imagination

The choice to conceive an educational path on operational strategies that develop to the establishment of mathematical knowledge intended, not only from the technical-operational point of view but also as complex and articulated activity of the human thought (beauty and ecstasy in the synergic effect between reality and fantasy) has been made taking into consideration careful considerations on ability of the mind to deduce conjectures, retract or verify them in a frame of reference not strictly attributable to rationality as form of mathematical expression (in general and of geometrical expression in particular) for excellence.

The knowledge comes mainly as element of mediation between two types of thinking: *convergent thinking*, founded on experience, and *divergent thinking* that, going beyond the boundaries of reality, is placed in the perceptive, intuitive, creative dimension.

Starting from studies on different functions of the brain by Roger Sperry (Morabito 2004), it is possible to assert what are the criteria that govern the evolutionary process of the brain and their utility in the context of knowledge on creative thought. Specifically, it is studied that the human brain has the structure of a complex modular system where different functions are encoded by the integration of different sets of specific brain structures.

The functionality of the brain is an expression of a double thinking that expresses itself through two complementary procedures: the logical-deductive (analytical-deductive) reasoning and the intuitive (synthetic-inductive) reasoning which correspond to functionally different procedures of the actions of the two cerebral hemispheres. It is due to this ability of the mind that humans changed the world in the century by creating tools helpful to their survival and evolution.

Roger Sperry showed that the correct correlation and the complementarity between these two modes of expression of the thinking had taken on great importance; in summary between rationality and creativity. To try alternative solutions, to think individual logical mathematical concepts, to create particular and innovative geometrical constructions mean to express an essential and constituent part of humans. So to think the mathematics trough these concepts means to carry out the data of reality filtering through always-new intuitions, which continually modify the essential characteristics of them.

## 2   Geometric figures, creativity and art

The main functions of drawing are to stimulate creativity and to construct individual imagines of reality. When we draw lines that create geometric shapes, for example, we are looking trough the eyes of the mind beyond the simple specified form, attributing a series of geometries, dynamic tensions, proportions to it.  In summary there is a close continuous dynamic interaction between what is graphically symbolized in the geometrical shape we see and what constantly changes and is imprinted as image in our mind. A synergic and dynamic connection between what is drawn on paper and what is formed and changes in our mind begins.

It is a special capacity of humans to see, in the complexity of shapes that surround them, the essential structures which constitute the shapes and interpret them in connection with their individual and subjective abilities. We see through our eyes, but we discern and know full well these images through our brain.  The mechanisms at the basis of this concept are complex, detailed and subject of constant studies.

The continuous interaction between idea and design allows, in summary, the materialization of the form in its essence.

In his Philebus, Plato wrote: "For when I say beauty of form, I am trying to express, not what most people would understand by the words, such as the beauty of animals or of paintings, but I mean, says the argument, the straight line and the circle and the plane and solid figures formed from these by turning-lathes and rulers and patterns of angles; perhaps you understand. For I assert that the beauty of these is not relative, like that of other things, but they are always absolutely beautiful by nature"; (Plato, Philebus 51c).

Now to a certain extent, this seems somewhat outdated but the conviction of profound necessity regard to the "control of the shape" through the correct use of the mental representation hold steadies.

## 3   From art to geometry

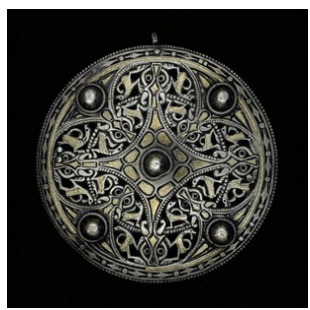

Fig. 1. Strickland Brooch - 9thC(mid) - http://www.britishmuseum.org

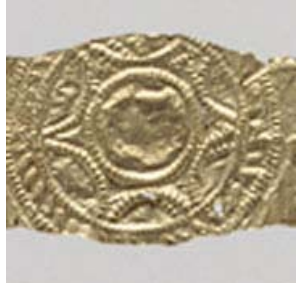

Fig. 2. Diadem, 400–450 Ostrogothic - The Metropolitan Museum of Art

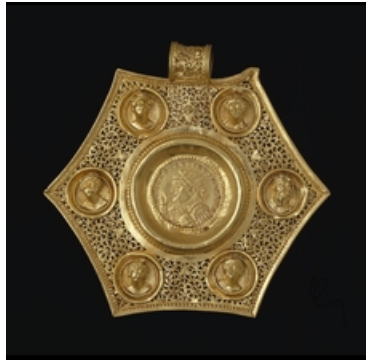

Fig. 4. Coin-set pendant - Late Antique, 4th century AD - http://www.britishmuseum.org

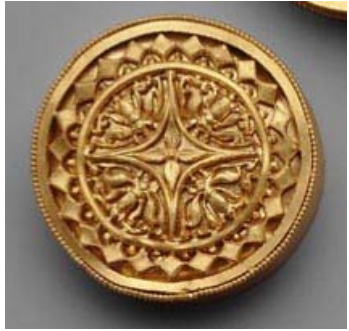

Fig. 5. Pair of Earplugs Ornamented - early 9th–14th century - The Metropolitan Museum of Art

The images above come from the collections of some of the most important museums in the world and depict historical and artistic artefacts from different periods. If we observe carefully the inner geometric shapes used for the decoration of the objects, they are all referable to the same type of curve now called *hypocycloid*.

According to the classical definition, the hypocycloids are the plane curves generated by a point on the circumference of a circle as it rolls without slipping along the inside of a circle (Fig. 6.).

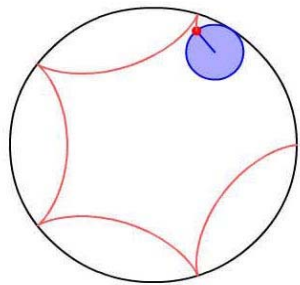

Fig. 6. Hypocycloid's construction

If the smaller circle has radius *r* and the larger circle has radius *R = kr*, where *k* is an integer, then the curve is closed and has *k* cusps.
During the centuries, famous mathematicians were interested on alternative methods for generate hypocycloids and, in particular, the "tricuspid hypocycloid" (curve of third class and fourth order[1]). According to the classical definition, if the radius *r* of the movable circumference is 1/3 (or, that is equivalent, 2/3) of the radius *R* of the fixed circumference, the curve formed is the tricuspid hypocycloid.

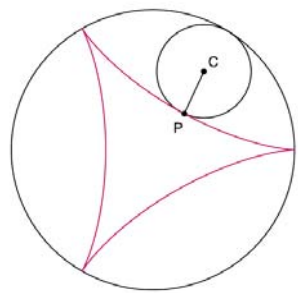

Fig. 7. The tricuspid hypocycloid

The parametric equations of this curve are:

$$\begin{cases} x = r(2\cos\alpha + \cos 2\alpha) \\ y = r(2\sin\alpha - \sin 2\alpha) \end{cases}. \qquad (1)$$

In geometry, given a triangle *ABC* and a point *P* on its circumcircle, the three closest points to *P* on lines *AB*, *AC*, and *BC* are collinear. The line through these points is the "Simson-Wallace line".[2]
Nowadays interactive geometry software allows us to create and then manipulate geometric constructions in plane geometry and so we are able to build and verify this property and the others that follow.

---

[1] The class of an algebraic curve is the number of its tangents that can be drawn from a generic point.
[2] In 1814 François-Joseph Servois (1767-1847) used the Simson-Wallace line which himself attributed to Robert Simson (1687-1768), however illustrating his assertion with the expression "je crois", in (Servois 1814). In 1822 in his paper *Traité des propriétés projectives des figure*, Jean-Victor Poncelet (1788-1867) reiterated that Servois had ascribed the theorem to Simson, but he omitted the expression "je crois"; (Poncelet 1822). In a note in (Servois 1814), Gergonne proved in detail this theorem using algebraic methods. In any case the Simson-Wallace line is present in Wallace 1798, v. II, p. 111 by William Wallace (1768-1843).

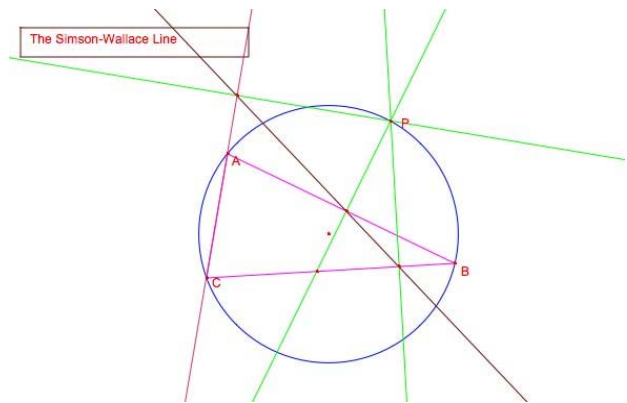

Fig. 8. The Simson-Wallace line

It is possible to use our dissertation as a starting point to conceive a geometric laboratory that, drawing inspiration from artistic works, get to create geometric shapes provided with fascinating symmetries.[3]

In 1857 the famous Swiss mathematician Jakob Steiner (1796-1863) showed that the envelope of the Simson Wallace line while the point moves on the circle is a tricuspid hypocycloid;[4] (J. Steiner 1857). If the given triangle is equilateral, the centre of the hypocycloid coincides with the centre of the given circle.

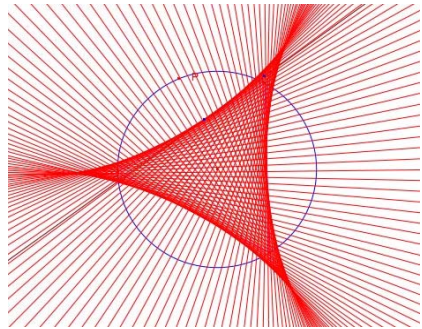

Fig. 9. The tricuspid hypocycloid as envelope of the Simson-Wallace line starting from an equilateral triangle

---

[3] The following pictures are created through Cabri II Plus; http://www.cabri.com/
[4] Actually the Italian mathematician Luigi Cremona (1830-1903) showed the equivalence of the curve created by Steiner with the tricusped hypocycloid; (Cremona 1864).

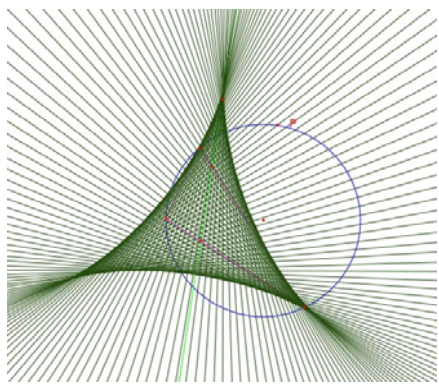

Fig. 10. The tricuspid hypocycloid as envelope of the Simson-Wallace line starting from a scalene triangle

In 1895, John Edward Aloysius Steggall (1855-1935), whose research interests were particularly in the geometry of the triangle, investigated the envelope of the Wallace lines of an inscribed polygon. He gave the equation of the Simson-Wallace line in the general case and some facts about regular polygons; he gave the condition that the envelope be a four-cusped hypocycloid; (Steggall 1895).

He started by building the Simson-Wallace line of a quadrangle: from a given point, perpendicular lines be let fall on the four Simson-Wallace lines formed from the four triangles made by taking every three of four points concyclic with the first; the feet of these perpendicular lines lie on a straight line that is called the Simson-Wallace line of the quadrangle starting from the four points. If we take an inscribed quadrangle, and omit each vertex in turn, we obtain four triangles; and the feet of the perpendiculars dropped from a point of the circle upon its four Simson-Wallace lines regard to these triangles lie on the Simson-Wallace line of the quadrangle. If the inscribed quadrangle is a square, then its Simson-Wallace line envelope a four-cusped hypocycloid.

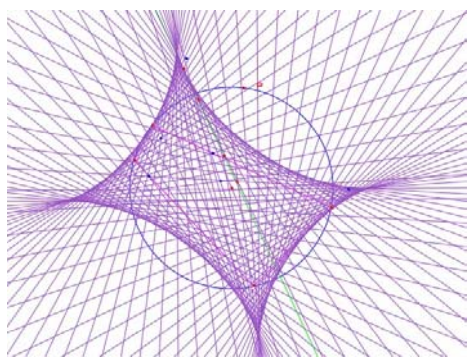

Fig. 11. The envelope of the Simson-Wallace line of an irregular quadrangle

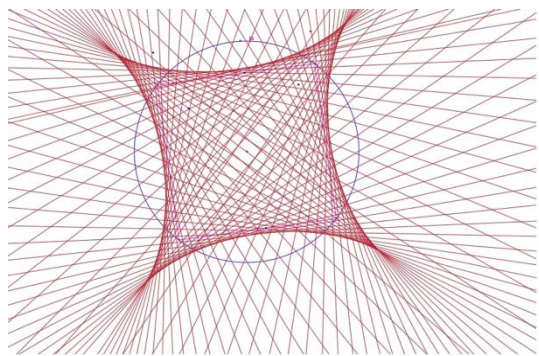

Fig. 12. The four-cusped hypocycloid (also called "astroid") as envelope of the Simson-Wallace line of a square

In 1919, David F. Barrow (1888-1970) came to the generalization that the Simson-Wallace line of an *n*-sided regular polygon envelops an *n*-cusped hypocycloid. After defining the Simson-Wallace line of a triangle, it is possible to define the Simson-Wallace line of an *n*-sided regular polygon: if we omit each vertex in turn, we obtain *n* irregular polygons with *n-1* sides. From a given point *P*, perpendicular lines be let fall on the *n*-1 Simson-Wallace lines formed from the *n* irregular *n*-1 sided polygons made by taking every *n*-1 of *n* points concyclic with the first; the feet of these perpendicular lines lie on a straight line that is called the Simson-Wallace line of the given *n*-sided polygon; (Barrow 1919).

If the inscribed *n*-sided polygon is regular, its Simson-Wallace line envelops an *n*-cusped hypocycloid.

It is possible to adopt similar procedure for both regular and irregular polygons, but if the given *n*-sided polygon is regular, the envelope created is provided with fascinating symmetries.

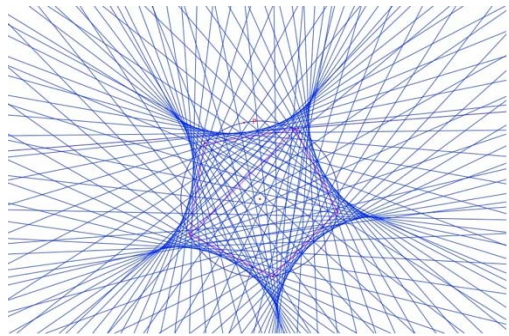

Fig. 13. The five-cusped hypocycloid as envelope of the Simson-Wallace line of a regular pentagon

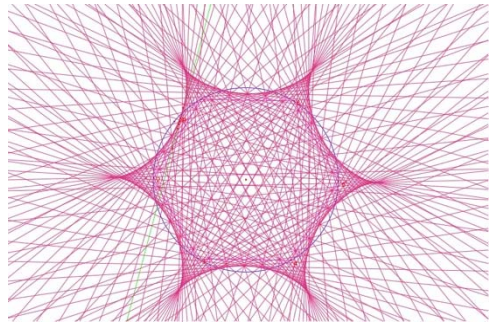

Fig. 14. The six-cusped hypocycloid as envelope of the Simson-Wallace line of a regular hexagon

## 4   An example of geometrical construction

We want to describe as example the construction of the envelope of the Simson-Wallace line starting from a square, through Cabrì. The other constructions in this paper are easily deduced from this.

| | |
|---|---|
| Construct a square and its circumcircle:<br>Lines -> Regular Polygon;<br>Curves ->Circle.<br>Take a point P on the circle:<br>Points ->Point on Object (P).<br>To give a name to the vertices of the square, use:<br>Text and Symbols -> Label.<br>Construct a triangle of vertices I, II, III (in green):<br>Lines -> Triangles.<br>Make the three sides longer:<br>Lines -> Line. | 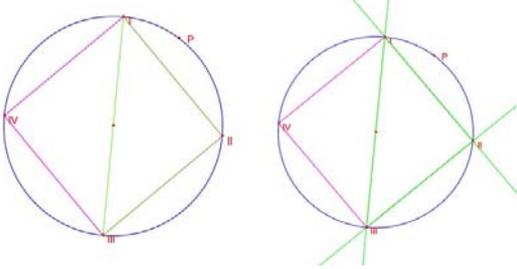 |
| Construct the Simson-Wallace line of the triangle with vertices I, II, III; from P draw the perpendicular lines on the three sides:<br>Constructions -> Perpendicular Line.<br>Take the feet of the perpendicular lines (in blue):<br>Points -> Intersection Point(s). | 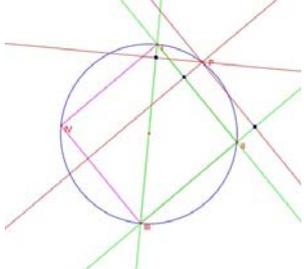 |
| To make more clear the construction, it is possible to hide some lines:<br>Attributes -> Hide/Show.<br>Draw the line passing trough the feet (in blu); this is the Simson-Wallace line:<br>Lines -> line. | 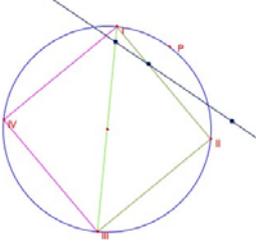 |
| Now we have to repeat the same construction for the triangle with vertices II, III, IV to obtain the second Simson-Wallace line. | 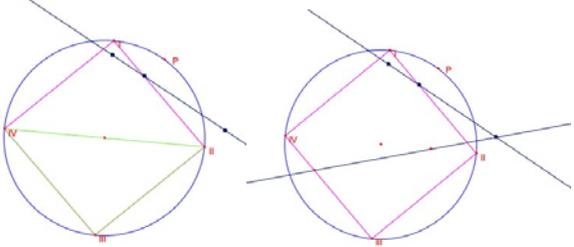 |
| By repeating the construction for the triangles with vertices III, IV, I and IV, I, II, we obtain the other two Simson-Wallace lines (in blue). | 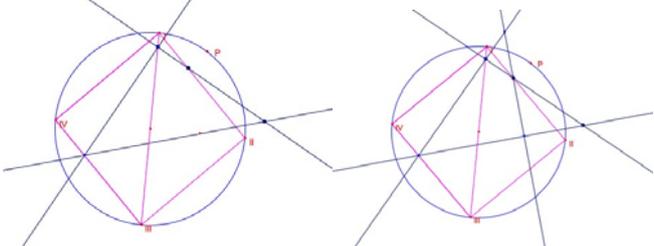 |
| Then we have to draw the perpendicular lines (in green) from P to the four Simson-Wallace lines and take the four feet A, B, C, D (in orange). | 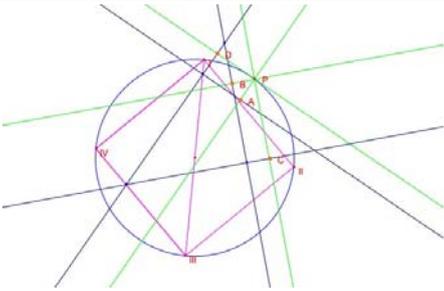 |
| Trough A, B, C, D passes a line that is the Simson-Wallace line of the given square (in orange).<br>To construct the envelope (as in Fig. 12), use:<br>Construction –Locus<br>setting Options -> Preferences -> Loci Options on 100 as number of objects in locus. | 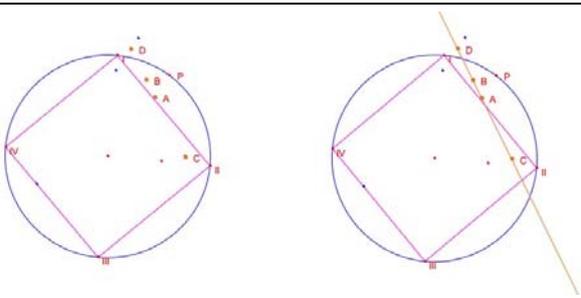 |

## 5  Conclusions

Over the years, the teaching of mathematics, in general, and geometry, in particular, were subject to many changes. For example it is possible to consider the whole diatribe of the 70s on modern mathematics and the dissemination of ideas of some members of the Bourbaki group; in particular the famous statement by Jean Dieudonné (1906-1992): "À bas Euclide"; (see Mammana, Villani 1998). Too often there is a distorted way of thinking mathematics (and geometry) attributing to it the label of sterile formalism.
De Finetti (De Finetti 1976), in this regard, said that unfortunately Mathematics had, and has, in the collective imagination the connotation of static and unchanging science; a greater attention to its historical development would help to change this misconception.
We think that combining in the learning process some historical notes, that highlight the continuous evolution of the mathematics, and suggestions aroused by the beauty of fascinating symmetries can stimulate the interest and the learning of simple or complex concepts, within the geometry.
So we want to conclude saying that appearance and contemplative have their importance in creating geometrical images of certain beauty.

## References

[1] BARROW, D. On the Envelope of the Wallace Lines of an Inscribed Quadrangle, The American Mathematical Monthly, vol. 26, n.3, 1919, pp.108-111.
[2] CREMONA, L. Sur l'hypocycloide à trois rebroussements. Journal für die reine und angewandte Mathematik, Issue 64.
[3] D'AMORE, B. Giocare con la matematica, Bologna, Archetipolibri, 2009.
[4] DE FINETTI, B. Dall'utopia all'alternativa, Milano, Franco Angeli Editore, 1976.
[5] GARDNER, M. I misteri della magia matematica. Firenze, Sansoni, 1985.
[6] GUILFORD, J. P. Creativity. The American Psychologist, 1950.
[7] MORABITO, C. La mente nel cervello. Roma-Bari, Laterza Editori, 2004.

[8] MAMMANA, C., VILLANI, V., Geometry and geometry-teaching through the ages, in Perspectives on the Theaching of Geometry for the 21st century, ICMI study, Kluwer 1998

[9] PLATO. Plato in Twelve Volumes, Vol. 9 translated by Harold N. Fowler. Cambridge, MA, Harvard University Press; London, William Heinemann Ltd. 1925.

[10] PONCELET, J.V. Traité des propriétés projectives des figures. Bachelier, Paris, 1822, p. 261.

[11] SERVOIS, F.J. Géométrie pratique. Problème. Prolonger une droite accessible audelà d'un obstacle qui borne la vue, en n'employant que l'équerre d'arpenteur, et sans faire aucun chaînage?. Annales de Mathématiques pures et appliquées, tome 4, 1813-1814, p. 250-253.

[12] STEINER, J. Über eine besondere Curve dritter Klasse (und vierten Grades). Journal für die reine und angewandte Mathematik, Berlin, p. 231-237, 1857.

[13] STEGGALL, J.E.A. On the envelope of the Simson line of a polygon, Proceedings of The Edinburgh Mathematical Society, Vol. 14 / February 1895, pp. 122-126, 1895.

[14] WALLACE, W. Mathematical lucubrations. Mathematical Repository, by T. Leybourn, printed by W. Glendinning, London, 1798.